\documentclass[a4paper,12pt,oneside]{article}
\usepackage[left=2.5cm,right=2.5cm,top=2.5cm,bottom=2.5cm]{geometry}
\usepackage[utf8]{inputenc}
\usepackage{t1enc}
\usepackage[T1]{fontenc}
\usepackage{xcolor}

\usepackage{amsmath}
\usepackage{amsfonts}
\usepackage{amssymb}
\usepackage{amsthm}
\usepackage{url}
\title{Lower bounds for coefficients of certain cyclotomic polynomials}
\author{Ákos Borsányi\\
Eötvös University, Pázmány Péter sétány 1/c, Budapest 1117, Hungary\\
E-mail: akos.borsanyi@gmail.com}
\date{2024.03.04.}

\begin{document}
	\maketitle
	

\renewcommand{\thefootnote}{}



\renewcommand{\thefootnote}{\arabic{footnote}}
\setcounter{footnote}{0}

\begin{abstract}
In this paper I prove a conjecture which gives a lower bound
for the largest absolute value of the coefficients of the $n$-th
cyclotomic polynomial for some $n$. Moreover this estimate is essentially sharp.
\end{abstract}

	\section{Introduction}
	
	The $n$-th cyclotomic polynomial is defined by the next product.
	\[\Phi_n(z)=\prod_{\substack{l=1\\ (l, n)=1}}^{n}\left(z-\varepsilon^l\right),\qquad\text{where }\varepsilon=e^{\frac{2i\pi}{n}}.\]
	It is well known that
	\begin{equation}\label{elso}
		\Phi_n(z)=\prod_{d\mid n} \left(z^\frac{n}{d}-1\right)^{\mu(d)}
	\end{equation}
	(Möbius's inversion formula).
	
	We denote by $A(n)$ the largest absolute value of the coefficients of $\Phi_n(z)$. From \eqref{elso} it is easy to see that if $p$ is a prime divisor of $n$ then $\Phi_{np}(z)=\Phi_n\left(z^p\right)$, so $A(np)=A(n)$. Thus, if we investigate $A(n)$ we may assume that $n$ is square free.
	In the sequel we assume that $2<p_1<p_2<\ldots<p_k$ are primes and $n=\prod_{l=1}^{k}p_l$. In \cite{bateman}
	it is proved that 
	\begin{equation}\label{masodik}
		A(n)\le\prod_{l=1}^{k-2}p_l^{2^{k-1-l}-1}
	\end{equation}
	moreover I have shown in my thesis \cite{szakdolgozat} that
	\begin{equation}\label{harmadik}
		A(n)\le c_k\prod_{l=1}^{k-2}p_l^{2^{k-1-l}-1}
	\end{equation}
	where $c_1=c_2=1$, $c_3=c_4=\frac{3}{4}$ and $c_k=\left(\frac{3}{8}\right)^{2^{k-5}}$ if $k\ge 5$.
	
	It is not obvious but not so difficult to see that
	\[\prod_{l=1}^{k-2}p_l^{2^{k-1-l}-1}\le n^{\frac{2^{k-1}}{k}-1}\]
	hence
	\begin{equation}\label{negyedik}
		A(n)\le c_k n^{\frac{2^{k-1}}{k}-1}.
	\end{equation}
	
	If the primes $p_l$ are closed to each other, that is $p_k-p_1$ is small then \eqref{negyedik} is not much weaker than \eqref{harmadik}. But note that by \eqref{masodik} the size of $A(n)$ is independent of the size of $p_{k-1}$ and $p_k$. Because of this if $p_1, \ldots, p_{k-2}$ are fixed and $p_{k-1}$ and $p_k$ tends to infinity we get $A(n)=O(1)$ and hence there do not exist any constants $d_k$ for which
	\begin{equation}\label{otodik}
		A(n)\ge d_k n^{\frac{2^{k-1}}{k}-1}
	\end{equation} 
	for every $n$.
	
	On the other hand it has been conjectured (see \cite{bateman}) that \eqref{otodik} holds for infinitely many $n$, where $\omega(n)=k$ and $d_k$ is fixed for every $k$. The purpose of this paper is to give a proof of this conjecture using the next recent result due to Maynard \cite{maynard}. For each $k$ positive integer there exists a bound $L_k$ for which there exist infinitely many $k$-tuples $p_1<p_2<\ldots<p_k\le p_1+L_k$ consisting of primes.
	
	\section{The proof of the conjecture}
	
	Let $k$ be fixed, $p_1<p_2<\ldots<p_k$  and $n=\prod_{l=1}^k p_l$. By \eqref{elso} for any real $x$ for which $x\neq \frac{a\pi}{n}$ where $(a,n)\neq 1$ we have 
	\[\left|\Phi_n\left(e^{2ix}\right)\right|=\prod_{d\mid n}\left|e^{2dix}-1\right|^{\mu\left(\frac{n}{d}\right)}\]
	and  for any real $x$ for which $x\neq \frac{ap_k\pi}{n}$ where $\left(a, \frac{n}{p_k}\right)\neq 1$ we have
	\begin{equation}\label{hatodik}
		\begin{split}
		\left|\Phi_n\left(e^{2ix}\right)\right|
		&=\left|\sum_{l=0}^{p_{k-1}}e^{2p_1\cdots p_{k-1}lix}\right|
		\cdot\prod_{\substack{d\mid n\\ d\neq n\\ d\neq n/p_k}}\left|e^{dix}-e^{-dix}\right|^{\mu\left(\frac{n}{d}\right)}=\\
		&=\left|\sum_{l=0}^{p_{k-1}}e^{2\frac{n}{p_k}lix}\right|
		\cdot\prod_{\substack{d\mid n\\ d\neq n, n/p_k}}\left|\sin dx\right|^{\mu\left(\frac{n}{d}\right)}	
		\end{split}
	\end{equation}
	
	Denote by $k_1$ the integer part of $\frac{k}{2}$ and set $x=\frac{a\pi}{p_1\cdots p_{k-1}}$ where for the first time $a$ is an arbitary integer relatively prime to $\prod_{l=1}^{k-1}p_l$. Let $\varepsilon=e^{\frac{2i\pi}{p_1\cdots p_{k-1}}}$. Now by \eqref{hatodik}
	\begin{multline}\label{hetedik}
		\left|\Phi_n(\varepsilon^a)\right|
		=p_k\prod_{\substack{d\mid n\\ d\neq n, n/p_k}}\left|\sin\frac{ad\pi}{p_1\cdots p_{k-1}}\right|^{\mu\left(\frac{n}{d}\right)}=\\
		=p_k\frac
		{\prod_{m=1}^{k_1}\prod_{1\le i_1<\ldots<i_{2m-1}\le k-1}\left|\sin\frac{a\pi}{\prod_{s=1}^{2m-1}p_{i_s}}\right|
		\cdot\prod_{m=1}^{k_1}\prod_{1\le i_1<\ldots<i_{2m}\le k-1}\left|\sin\frac{ap_k\pi}{\prod_{s=1}^{2m}p_{i_s}}\right|}
		{\prod_{m=1}^{k_1}\prod_{1\le i_1<\ldots<i_{2m}\le k-1}\left|\sin\frac{a\pi}{\prod_{s=1}^{2m}p_{i_s}}\right|
		\cdot\prod_{m=1}^{k_1}\prod_{1\le i_1<\ldots<i_{2m-1}\le k-1}\left|\sin\frac{ap_k\pi}{\prod_{s=1}^{2m-1}p_{i_s}}\right|}.
	\end{multline}

	Note that in case when $k=2k_1$ the $\prod_{1\le i_1<\ldots<i_{k}\le k-1}\left|\sin\frac{ap_k\pi}{\prod_{s=1}^{k}p_{i_s}}\right|$ and the\linebreak $\prod_{1\le i_1<\ldots<i_{k}\le k-1}\left|\sin\frac{a\pi}{\prod_{s=1}^{k}p_{i_s}}\right|$ products are empty.
	
	According to Maynard's theorem mentioned in the introduction there exists $L_k=\colon L$ such that 
    \begin{equation}\label{nyolcadik}p_1<\ldots<p_k\le p_1+L_k=p_1+L 
    \end{equation}for infinitely many ($p_1,\ldots, p_k$) $k$-tuples where $p_1, \ldots, p_k$ are all primes. 
	
	In the sequel we suppose that the primes $p_l$ ($1\le l\le k$) and $L$ are satisfying (\ref{nyolcadik}). For $1\le l\le k$ we define the integers $j_l$ by the equality $p_l=p_1+2j_l$. Now we distinguish the next two cases.
	
	\textbf{Case 1.} $j_2$ is odd i.e. $4\mid p_1+p_2$.
	
	\textbf{Case 2.} $j_2$ is even i.e. $4\mid p_1+p_2-2$.
	
	We carry out the proof in case 1. 
	
	Let $f_l(p_1)=\frac{p_1^{2l-1}+p_2^{2l-1}}{4}=\frac{p_1^{2l-1}+(p_1+2j_2)^{2l-1}}{4}$ ($1\le l\le k_1$) and define $a$ by $a=\sum_{l=1}^{k_1} f_l(p_1)$. The essence of the proof is to investigate the factors of the fraction \eqref{hetedik}. The main observation is that if $a, b$ and $M$ are integers, $M>0$ and $a\equiv b (\mod M)$ then $\left|\sin\frac{a\pi}{M}\right|=\left|\sin\frac{b\pi}{M}\right|$.
	
	First of all we shall show that if $m$ and $i_1,i_2,\ldots\,i_{2m-1}$ are fixed and $p_1$ tends to infinity then $\left|\sin\frac{a\pi}{\prod_{s=1}^{2m-1}p_{i_s}}\right|\to 1$. It is easy to see that if $l<m$ then $\frac{p_1^{2l-1}+p_2^{2l-1}}{4\prod_{s=1}^{2m-1}p_{i_s}}$ tends to $0$ and $\lim_{p_1\to\infty}\frac{p_1^{2m-1}+p_2^{2m-1}}{4\prod_{s=1}^{2m-1}p_{i_s}}=\frac{1}{2}$. Now let $l\ge m+1$. Then
 \begin{align*}
	f_l(p_1)=\frac{p_1^{2l-1}+p_2^{2l-1}}{4}=\frac{p_1^{2l-2}+(p_1+2j_2)^{2l-1}}{4}=\\
 =\frac{p_1+(2l-1)j_2}{2}p_1^{2l-1}+\sum_{t=2}^{2l-1}\binom{2l-1}{t}2^{t-2}j_2^t p_1^{2l-1-t}
 \end{align*}
	that is
	\begin{equation}\label{kilencedik}
		f_l(p_1)=\frac{p_1+(2l-1)j_2}{2}p_1^{2l-2}+\sum_{t=2}^{2l-3}\binom{2l-1}{t}2^{2l-3-t}j_2^{2l-1-t}p_1^t.
	\end{equation}
	
	Now
	\begin{equation}\label{tizedik}
		\frac{p_1+(2l-1)j_2}{2}p_1^{2l-2}\equiv \frac{p_1+(2l-1)j_2}{2}\left(p_1^{2l-2}-p_1^{2l-2m-1}\prod_{s=1}^{2m-1}(p_1+2j_{i_s})\right)\mod\prod_{s=1}^{2m-1}p_{i_s}.
	\end{equation}
	
	Observe that the right hand side of the congruence \eqref{tizedik} is a polynomial of $p_1$ with intger coefficients. So, because of \eqref{kilencedik} and \eqref{tizedik} we get $f_l(p_1)\equiv g_l(p_1)\mod\prod_{s=1}^{2m-1}p_{i_s}$ where $g_l$ is a polynomial of $p_1$ with integer coefficients and of degree at most $2l-2$.
	Since $g_l(x)\in\mathbb{Z}[x]$ and furthermore $\prod_{s=1}^{2m-1}p_{i_s}=\prod_{s=1}^{2m-1}(p_1+2j_{i_s})$ is a monic polynomial of $p_1$ with also integer coefficients $g_l(p_1)\equiv h_l(p_1)\mod \prod_{s=1}^{2m-1} p_{i_s}$ where $\deg h_l\le 2m-2$ and the magnitude of the coefficients of $h_l$ depends only on $k$. Therefore
	\begin{equation}\label{tizenegyedik}
		\left|\sin\frac{a\pi}{\prod_{s=1}^{2m-1}p_{i_s}}\right|
		=\left|\sin\frac{\sum_{l=1}^{k_1}f_l(p_1)\pi}{\prod_{s=1}^{2m-1}p_{i_s}}\right|
		=\left|\sin\frac{\sum_{l=1}^{m-1}f_l(p_1)\pi+f_m(p_1)\pi+\sum_{l=m+1}^{k_1}h_l(p_1)\pi}{\prod_{s=1}^{2m-1}p_{i_s}}\right|
	\end{equation}
	Since the degree of the polynomial $\sum_{l=1}^{m-1}f_l(p_1)+\sum_{l=m+1}^{k_1}h_l(p_1)$ is at most $2m-2$,
	\[\lim_{p_1\to\infty}\frac{\sum_{l=1}^{m-1}f_l(p_1) + f_m(p_1) + \sum_{l=m+1}^{k_1}h_l(p_1)} {\prod_{s=1}^{2m-1}p_{i_s}}
	=\lim_{p_1\to\infty}\frac{f_m(p_1)}{\prod_{s=1}^{2m-1}p_{i_s}}=\frac{1}{2},\]
	so because of \eqref{tizenegyedik} $\lim_{p_1\to\infty}\left|\sin\frac{a\pi}{\prod_{s=1}^{2m-1}p_{i_s}}\right|=1$.
	
	Now consider $\left|\sin\frac{ap_k\pi}{\prod_{s=1}^{2m}p_{i_s}}\right|$ where $1\le m\le k_1-1$ if $k=2k_1$ and $1\le m\le k_1$ if $k=2k_1+1$ and $m$ and $i_1,\ldots, i_{2m}$ are fixed. Similarly to the previous argument if $l\ge m+1$ then $f_l(p_1)\equiv g_l(p_1)\mod\prod_{s=1}^{2m}p_{i_s}$ where $g_l(x)\in\mathbb{Z}[x]$ and $\deg g_l=2l-2$. So, because $\prod_{s=1}^{2m}p_{i_s}=\prod_{s=1}^{2m}(p_1+2j_{i_s})$ is a monic polynomial of $p_1$ having integer coefficients, $p_kf_l(p_1)=(p_1+2j_k)f_l(p_1)\equiv(p_1+2j_k)g_l(p_1)\equiv h_l(p_1)\mod\prod_{s=1}^{2m}p_{i_s}$ where $h_l(x)\in\mathbb{Z}[x]$, $\deg h_l\le 2m-1$ and the magnitude of the coefficients of $h_l(x)$ depends only on $k$. Hence as previously
	\[\lim_{p_1\to\infty}\left|\sin\frac{ap_k\pi}{\prod_{s=1}^{2m}p_{i_s}}\right|
	=\lim_{p_1\to\infty}\left|\sin\frac{p_k\sum_{l=1}^{m-1}f_l(p_1)+p_kf_m(p_1)+\sum_{l=m+1}^{k_1}h_l(p_1)}{\prod_{s=1}^{2m}p_{i_s}}\pi\right|=1.\]
	Consider a factor of the denominator of the fraction \eqref{hetedik} being of the form $\left|\sin\frac{a\pi}{\prod_{s=1}^{2m}p_{i_s}}\right|$. Since $a=\sum_{l=1}^{k_1}\frac{p_1^{2l-1}+(p_1+2j_1)^{2l-1}}{4}=\sum_{l=1}^{k_1}f_l(p_1)$ can be written as 
	$a=\sum_{l=1}^{m-1} f_l(p_1) + f_m(p_1) + \sum_{l=m+1}^{k_1} f_l(p_1)$ and as we have seen above that 
	$\sum_{l=m+1}^{k_1} f_l(p_1)\equiv g(p_1)\mod\prod_{s=1}^{2m}p_{i_s}$ where $g(x)\in\mathbb{Z}[x]$, $\deg g\le 2m-1$, it can be seen that $a\equiv h(p_1) \mod\prod_{s=1}^{2m}p_{i_s}$, where $h$ is a polynomial having degree at most $2m-1$ and the size of the coefficients of $h$ depends only on $k$. But taking the term $f_m(p_1)$ into consideration the leading coefficient of $h$ is of the form $b+\frac{1}{2}$ where $b\in\mathbb{Z}$ so the degree of $h$ is in fact $2m-1$ (as $b+\frac{1}{2}$ cannot be equal to $0$).
	
	Now because of this $\left|\sin\frac{a\pi}{\prod_{s=1}^{2m}p_{i_s}}\right|=\left|\sin\frac{h(p_1)\pi}{\prod_{s=1}^{2m}p_{i_s}}\right|=O\left(\frac{1}{p_1}\right)$ but not equal to $O\left(\frac{1}{p_1^2}\right)$.	Here the constant in the symbol $O$ depends only on $k$.
	
	Finally, investigate the factor $\left|\sin\frac{ap_k\pi}{\prod_{s=1}^{2m-1}p_{i_s}}\right|$ where $1\le m\le k_1$ and $1\le i_1<\ldots<i_{2m-1}\le k-1$ are fixed. We take apart $a$ in the following way.
	\begin{equation}\label{tizenkettedik}
		a=\sum_{l=1}^{m-2}f_l(p_1) + f_{m-1}(p_1)+\sum_{l=m}^{k_1}f_l(p_1)
	\end{equation}
	
	By \eqref{kilencedik} the degree of $p_k\sum_{l=1}^{m-2}f_l(p_1)=(p_1+2j_k)\sum_{l=1}^{m-2}f_l(p_1)$ is $2m-4$ and $\deg(p_1+2j_k)f_{m-1}(p_1)=2m-2$. Furthermore the leading coefficient of $(p_1+2j_k)f_{m-1}(p_1)$ is of the form $b+\frac{1}{2}$ where $b\in\mathbb{Z}$. Now consider the term $\frac{p_1+(2l-1)j_2}{2}p_1^{2l-2}(p_1+2j_k)$ of $f_l(p_1)(p_1+2j_k)$ where $m\le l\le k_1$.
	\begin{multline*}
		\frac{p_1+(2l-1)j_2}{2}p_1^{2l-2}(p_1+2j_k)=(j_kp_1+(2l-1)j_2j_k)p_1^{2l-2}+\frac{p_1+(2l-1)j_2}{2}p_1^{2l-1}\equiv\\
		\equiv j_kp_1^{2l-1}+(2l-1)j_2j_kp_1^{2l-2}-\frac{p_1+(2l-1)j_2}{2}\left(p_1^{2(l-m)}\prod_{s=1}^{2m-1}(p_1+2j_{i_s})-p_1^{2l-1}\right)\mod\prod_{s=1}^{2m-1}p_{i_s}
	\end{multline*}
	
	Denote this polynomial of $p_1$ by $g_l(p_1)$. It is easy to see that $g_l(x)\in\mathbb{Z}[x]$.
	
Since $\prod_{s=1}^{2m-1}(p_1+2j_{i_s})$ is a monic polynomial of $p_1$, by \eqref{kilencedik}
	\[(p_1+2j_k)f_l(p_1)\equiv g_l(p_1)+(p_1+2j_k)\sum_{t=0}^{2l-3}\binom{2l-1}{t} 2^{2l-3-t}j_2^{2l-1-t}p_1^t\equiv h_l(p_1)\mod\prod_{s=1}^{2m-1}p_{i_s}\] where $h_l\in\mathbb{Z}[x]$ and $\deg h_l\le 2m-2$.
	
	According to \eqref{tizenkettedik} and the previous arguments $ap_k=a(p_1+2j_k)\equiv h(p_1)\mod\prod_{s=1}^{2m-1}p_{i_s}$ and $\deg h\le 2m-2$ but the leading coefficient of $h$ is of the form $b+\frac{1}{2}$ ($b\in\mathbb{Z}$) so $\deg h=2m-2$. So we get
	\[\left|\sin\frac{ap_k\pi}{\prod_{s=1}^{2m-1}p_{i_s}}\right|=\left|\sin\frac{h(p_1)\pi}{\prod_{s=1}^{2m-1}p_{i_s}}\right|=O\left(\frac{1}{p_1}\right)\]
	where the constant in the symbol $O$ depends only on $k$. Because of the previous remark $\left|\sin\frac{ap_k\pi}{\prod_{s=1}^{2m-1}p_{i_s}} \right|\neq O\left(\frac{1}{p_1^2}\right)$. 
	
	it is not so difficult to see that the number of the factors in the denominator of \eqref{hetedik} is $2^{k-1}-1$, so by \eqref{hetedik} $|\Phi_n(\varepsilon^a)|>c_kp_1^{2^{k-1}}$ for some $c_k$.
	
	Finally because $nA(n)\ge (\varphi(n)+1)A(n)\ge |\Phi_n(\varepsilon^a)|> c_kp_1^{2^{k-1}}$,
	\[A(n)>d_k n^{\frac{2^{k-1}}{k}-1}\]
	for some $d_k$ and this completes the proof of the conjecture in case 1.
	
	The beginning of the proof in case 2.
	For $1\le l\le k_1$ let 
	\[f_l(p_1)=\frac{p_1^{2l-1}+(p_2-2)^{2l-1}}{4}=\frac{p_1^{2l-1}+(p_1+2(j_2-1))^{2l-1}}{4}\]
	and $a=\sum_{l=1}^{k_1}f_l(p_1)$. Now
	\begin{multline*}
		f_l(p_1)=\frac{p_1+(2l-1)(j_2-1)}{2}p_1^{2l-2}+\sum_{t=2}^{2l-1}\binom{2l-1}{t}2^{t-2}(j_2-1)^t p_1^{2l-1-t}=\\
		=\frac{p_1+(2l-1)(j_2-1)}{2}p_1^{2l-2}+\sum_{t=0}^{2l-3}\binom{2l-1}{t}2^{2l-3-t}(j_2-1)^{2l-1-t} p_1^t
	\end{multline*}
and from this the proof of the conjecture is the same as in case 1.

\end{document}